\documentclass[letterpaper,journal]{IEEEtran}

\usepackage{amsmath,amssymb,amsfonts}
\usepackage{graphicx}
\usepackage{booktabs}
\usepackage{algorithm}
\usepackage{algpseudocode}
\usepackage{tikz}
\usetikzlibrary{arrows.meta,positioning,calc,fit,backgrounds}
\usepackage{cite}
\usepackage{bm}

\newcommand{\R}{\mathbb{R}}
\newcommand{\norm}[1]{\left\lVert#1\right\rVert}

\begin{document}
\raggedbottom

\title{OCP-Based Nonlinear Model Predictive Control for Quadrotor Trajectory Tracking}

\author{Shichao~Lv, Chuanzhi~Lv, Hongxia~Wang, and Huanshui~Zhang%
\thanks{This work was supported by the International (Regional) Cooperation and Exchange Program of the National Natural Science Foundation of China (W2611094), the Original Exploratory Program Project of National Natural Science Foundation of China (62450004), the Joint Funds of the National Natural Science Foundation of China (U23A20325), and Key R\&D Program of Shandong Province of China (2025CXPT099). (Corresponding author: Huanshui Zhang).}%
\thanks{Shichao Lv, Chuanzhi Lv, Hongxia Wang, and Huanshui Zhang are with the College of Electrical Engineering and Automation, Shandong University of Science and Technology, Qingdao 266590, China (e-mail: Shichao Lv, lv\_sc2020@sdust.edu.cn; Chuanzhi Lv, lczbyq@sdust.edu.cn; Hongxia Wang, whx1123@126.com; Huanshui Zhang, hszhang@sdu.edu.cn).}%
\thanks{Huanshui Zhang is also with the School of Control Science and Engineering, Shandong University, Jinan, Shandong 250061, China.}%
\thanks{All authors are also with the Shandong Provincial Key Laboratory of Intelligent Sensing and Autonomous Control.}}

\markboth{PREPRINT}%
{OCP-Based Nonlinear Model Predictive Control for Quadrotor Trajectory Tracking}

\maketitle

\begin{abstract}
Fast online optimization is essential for nonlinear model predictive control (NMPC) on resource-limited aerial platforms. This paper presents an analytical-derivative NMPC solver for quadrotor trajectory tracking. For the state-affine translational model with held thrust and attitude commands, fourth-order Runge--Kutta integration yields an exact discrete-time prediction map. Forward--backward difference equations provide the gradient and exact Hessian of the finite-horizon cost condensed over the input sequence. These derivatives drive a nested optimal-control-based direction recursion, while a shifted warm start reuses the preceding solution. The resulting solver addresses the unconstrained tracking problem without a general-purpose nonlinear programming package. Gazebo experiments show comparable tracking accuracy to iLQR and two acados SQP configurations, with mean solver-time speedups of $1.97$--$2.30$. In onboard flight tests, the mean solver time is $0.3121$\,ms, corresponding to speedups of $3.37$--$5.77$ over the same baselines. The largest measured solver time is $3.7754$\,ms within a $20$-ms control period. These measurements cover the numerical solver and exclude estimation, communication, and command publication.
\end{abstract}


\section{Introduction}
\IEEEPARstart{U}{nmanned} aerial vehicles (UAVs) are used for inspection,
mapping, and aerial manipulation, where accurate tracking of three-dimensional
trajectories is often required~\cite{mellinger2011}. Their underactuated and
nonlinear dynamics make nonlinear model predictive control (NMPC) an attractive
tracking framework~\cite{mayne2014,nmpc_marine}. By solving a finite-horizon
optimal control problem at each sampling instant and applying its first input,
NMPC accounts for predicted tracking performance as the trajectory evolves.
Its effectiveness has been demonstrated in agile flight and
disturbance-rejection tracking%
~\cite{sun2022comp,recalde2025dq,gu2025taes,
andrien2024tac,sarvaiya2025hpa}.

Onboard NMPC typically operates under tight time limits for each
finite-horizon solve~\cite{cimen2004}. Continuation/generalized minimal residual methods reduce the work per
sample by updating the
Karush--Kuhn--Tucker conditions along the sampling sequence%
~\cite{wang2021icgmres}. Short stabilizing horizons,
state-dependent-coefficient reformulations, and code generation for
microcontrollers reduce the problem size or replace the nonlinear program with
simpler subproblems~\cite{gomaa2023cstb,omidi2025sdc,nguyen2025nmpcm}.
Analytical Koopman lifting leads to a quadratic program%
~\cite{rajkumar2025lra}, whereas event-triggered NMPC avoids solving the
optimization problem at every sampling instant~\cite{bao2026tie}. Although these approaches substantially improve real-time feasibility, their computational savings are generally obtained through problem reformulation,
model approximation, or fewer optimization updates. 

A complementary route to real-time NMPC is to reduce derivative and
optimization costs while retaining the nonlinear prediction model. The
finite-horizon NMPC problem solved at each sampling instant can be viewed in two
ways. The static viewpoint casts it as a nonlinear program and uses methods
such as interior-point, quasi-Newton, or sequential quadratic programming (SQP), with derivatives obtained
numerically or by automatic differentiation%
~\cite{nocedal,casadi,ipopt,lindqvist2020ra}. This is flexible but can
introduce software and computational overhead. The dynamic viewpoint derives
forward and backward difference equations
(FBDEs) from Pontryagin's maximum principle and computes derivatives through
state and costate recursions. It avoids a general-purpose nonlinear programming
package, although conventional first-order updates can converge slowly%
~\cite{chernousko1982}.

A recent line of work unifies the two viewpoints. Building on the optimization
framework rooted in optimal control of~\cite{zhang2024scis}, which develops a
superlinearly convergent algorithm from an optimal-control reformulation, the
OCP methods in~\cite{lv2024agv,lv2025constr} compute the cost gradient by
solving FBDEs, and the subsequent work~\cite{lv2026tie} introduces a virtual
optimal-control problem to derive a second set of FBDEs for Hessian evaluation.
The latter delivers a complete second-order solver that does not rely on
optimization packages and was validated on a linear--quadratic regulator and a
planar automatic guided vehicle, where it matched the accuracy of mature
solvers while reducing computation time by roughly an order of magnitude. That
solver, however, retains the first-order forward-Euler system model, because
higher-order discretization of general nonlinear dynamics makes the associated
derivative recursions substantially more involved. The resulting $O(h)$ global
integration error may therefore limit the attainable prediction and tracking
accuracy.

The quadrotor dynamics are affine in the state. With the input held
constant over a sampling interval, fourth-order Runge--Kutta (RK4) gives the
exact discrete-time model while preserving FBDE recursions. The resulting gradient and exact Hessian are used by the OCP
nested direction recursion to update the input sequence. The contributions
are as follows.
\begin{enumerate}
\item For real-time UAV trajectory tracking, an RK4-discretized model is
combined with FBDE recursions to obtain the analytical gradient and exact
Hessian of the finite-horizon cost with respect to the input sequence. The OCP
nested direction recursion uses these derivatives without a general-purpose
optimization package. 

\item A shifted warm start uses the preceding input sequence to initialize the optimization at each sampling instant. The applied input is discarded, the remaining inputs move forward
one stage, and the terminal reference input fills the last stage. This reduces
online iterations and solution time.

\item Gazebo and real-flight tests quantify the computational performance. In
the Gazebo comparison against the iterative
linear quadratic regulator (iLQR), acados SQP with the high-performance
interior-point method (acados SQP+HPIPM), and acados SQP with the online active-set QP solver (acados SQP+qpOASES),
Algorithm~\ref{alg:1} achieves a mean solution time of $0.1197$\,ms,
which is $1.97$--$2.30$ times faster than the baseline methods at comparable
tracking accuracy and with no recorded failures. The corresponding real-flight
mean is $0.3121$\,ms, giving speedup factors of $3.37$--$5.77$. The largest
measured time is $3.7754$\,ms, leaving sufficient solver-time margin for the
$20$-ms control period.
\end{enumerate}

\begin{figure*}[!t]
\centering
\resizebox{\textwidth}{!}{%
\begin{tikzpicture}[
  x=1mm,y=1mm,
  font=\scriptsize,
  >=Latex,
  flow/.style={->,line width=0.65pt,line cap=round,line join=round},
  crossing/.style={preaction={draw=white,line width=2.4pt}},
  block/.style={draw,rounded corners=1.5pt,align=center,
                inner sep=2.1pt,line width=0.7pt},
  refblock/.style={block,minimum height=12mm,fill=green!8,
                  draw=green!45!black,text width=20mm},
  stageone/.style={block,minimum height=17mm,fill=cyan!10,
                  draw=cyan!60!black,text width=26mm,line width=0.9pt},
  stagetwo/.style={block,minimum height=17mm,fill=orange!10,
                   draw=orange!70!black,text width=27mm,line width=0.9pt},
  stagethree/.style={block,minimum height=17mm,fill=orange!16,
                     draw=orange!80!black,text width=34mm,line width=0.9pt},
  warmblock/.style={block,minimum height=7mm,fill=yellow!10,
                    draw=orange!55!black,text width=42mm},
  updateblock/.style={block,minimum height=7mm,fill=violet!10,
                     draw=violet!65!black,text width=24mm,font=\tiny},
  innerblock/.style={block,minimum height=12mm,fill=blue!8,
                    draw=blue!55!black,text width=23mm},
  plantblock/.style={block,fill=green!9,draw=green!45!black,text width=22mm},
  imublock/.style={block,minimum height=11mm,fill=blue!10,
                  draw=blue!60!black,text width=10mm},
  lidarblock/.style={block,minimum height=11mm,fill=green!10,
                    draw=green!50!black,text width=14.5mm,font=\tiny},
  fastlioblock/.style={block,minimum height=8mm,fill=green!14,
                     draw=green!60!black,text width=14mm},
  fusionblock/.style={block,minimum height=8mm,fill=orange!15,
                     draw=orange!75!black,text width=22mm},
  methodgroup/.style={draw=cyan!60!black,fill=cyan!4,rounded corners=2.5pt,
                     line width=1.15pt,inner sep=3mm},
  group/.style={draw,dashed,rounded corners=2pt,line width=0.8pt,inner sep=3.5mm},
  signal/.style={font=\tiny,inner sep=0.6pt,outer sep=0pt}
]

\node[refblock] (trajectory) at (0,-0.5) {Trajectory\\reference\\$X^r,U^r$};

\node[stageone] (rk4) at (32,4.5)
  {\textbf{RK4 model}\\exact discretization\\forward state prediction};
\node[stagetwo] (fbde) at (67,4.5)
{\textbf{Analytical Derivatives}\\
FBDE-Based Gradient and Hessian};
\node[stagethree] (ocp) at (104,4.5)
  {\textbf{OCP Iterative Update}\\
  $d_0=(\mathcal R+H_i)^{-1}\nabla J_i;$\\
  $d_j=(\mathcal R+H_i)^{-1}(\nabla J_i$\\$+\mathcal R d_{j-1})$};
\node[warmblock] (warm) at (57,-10)
  {Shift $\mathcal U_{t-1}^{\star}$ and append $U^r_{N-1|t}$ to form
   $\mathcal U_t^0$};
\node[updateblock] (update) at (104,-10)
  {$\mathcal U_{i+1}=\mathcal U_i-d_i(\mathcal U_i)$};

\coordinate (refrouteL) at (16.5,16.5);
\coordinate (refrouteR) at (67,16.5);
\draw[flow] (trajectory.east) -- (16.5,-0.5) -- (refrouteL) --
  node[signal,below=0.9mm] {$X^r,U^r$} (refrouteR) -- (fbde.north);
\draw[flow] (rk4) -- node[signal,align=center,above=1.1mm]
  {$X_{0:N}$\\$\mathcal U_i$} (fbde);
\draw[flow] (fbde) -- node[signal,align=center,above=1.0mm]
  {$\nabla J_i$,\\$H_i$} (ocp);
\coordinate (ekftarget) at ($(rk4.south west)!0.18!(rk4.south east)$);
\coordinate (warmtarget) at ($(rk4.south west)!0.35!(rk4.south east)$);
\coordinate (looptarget) at ($(rk4.south west)!0.52!(rk4.south east)$);
\coordinate (loopbottomR) at (104,-15.5);
\coordinate (loopbottomL) at (looptarget |- loopbottomR);
\draw[flow] (ocp.south) -- (update.north);
\draw[flow,dashed] (update.west) --
  node[signal,below=0.8mm] {next sample} (warm.east);
\draw[flow] (update.south) -- (loopbottomR) -- (loopbottomL) -- (looptarget);
\coordinate (warminput) at (warmtarget |- warm.west);
\coordinate (warmcross) at (looptarget |- warm.west);
\draw[flow] (warm.west) -- ($(warmcross)+(1.1mm,0)$)
  .. controls ($(warmcross)+(0.75mm,1.25mm)$)
           and ($(warmcross)+(-0.75mm,1.25mm)$) ..
     ($(warmcross)+(-1.1mm,0)$)
  -- (warminput) -- (warmtarget);

\begin{scope}[on background layer]
\node[methodgroup,fit=(rk4)(fbde)(ocp)(warm)(update)(loopbottomR)(loopbottomL)
      (refrouteR),
      label={[font=\bfseries\scriptsize,cyan!55!black]above:
      Proposed Algorithm~\ref{alg:1}}] (outer) {};
\end{scope}

\node[innerblock] (attref) at (144,4.5)
  {Thrust--attitude\\command interface\\$T,\phi_d,\theta_d$};
\node[innerblock] (attctrl) at (144,-10.5)
  {Inner-loop\\attitude controller\\and motor allocation};
\coordinate (updatedoutput) at (127,-10);
\coordinate (outputbend) at (127,4.5);
\draw[flow] (update.east) --
  node[signal,above=0.9mm,pos=0.55] {$U_{0|t}^\star$}
  (updatedoutput) -- (outputbend) -- (attref.west);
\draw[flow] (attref) -- (attctrl);
\node[group,draw=blue!60!black,fit=(attref)(attctrl),
      label={[font=\bfseries\scriptsize,blue!60!black]above:Attitude loop}] (inner) {};

\node[plantblock,minimum height=17mm] (plant) at (186,9)
  {\includegraphics[width=17mm]{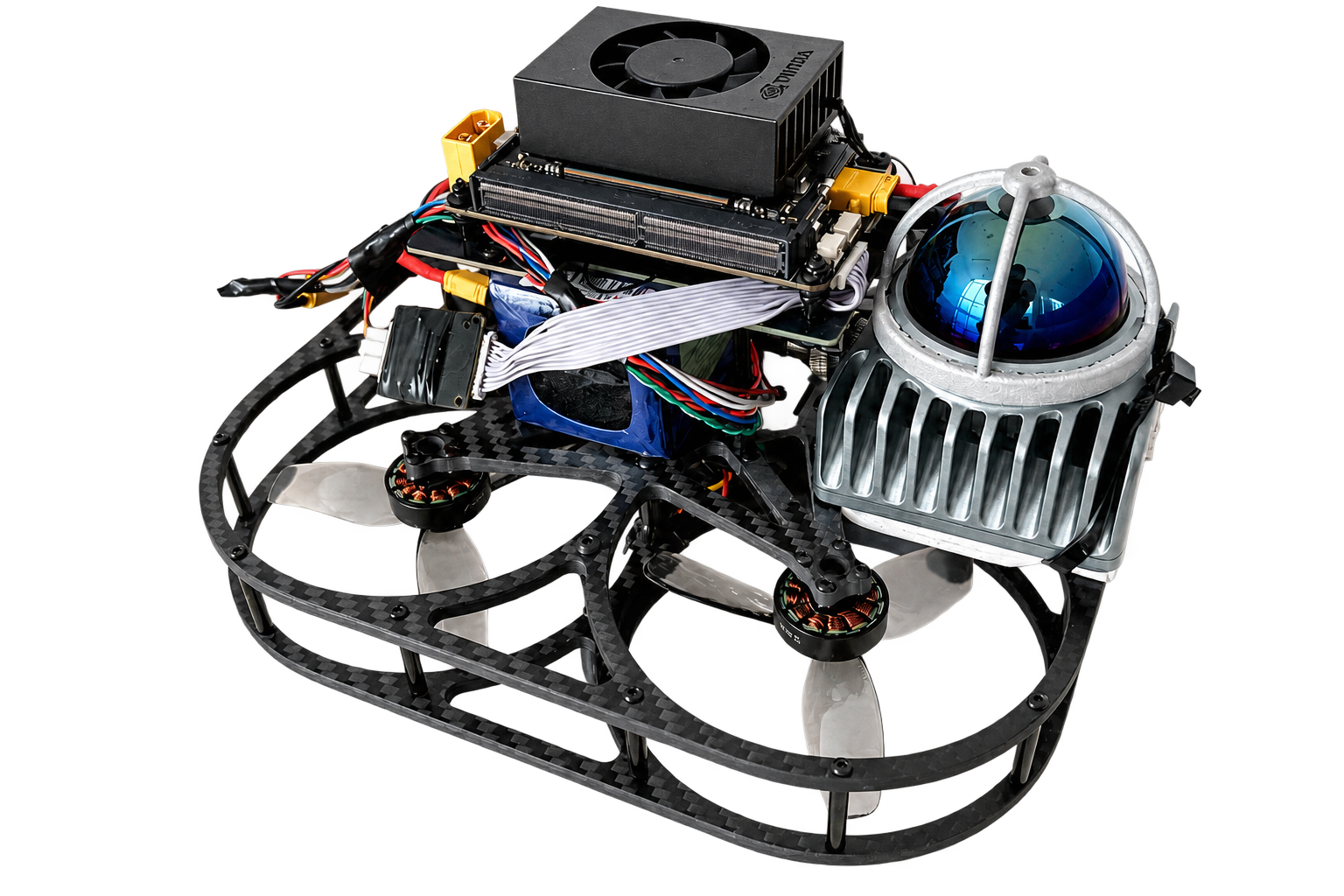}\\[-0.7mm]
   Quadrotor dynamics};
\node[imublock] (imu) at (180,-9) {IMU};
\node[lidarblock] (fastlio) at (193,-9) {LiDAR--IMU\\odometry\\(FAST-LIO)};
\node[fusionblock] (ekf) at (186,-21.5) {EKF state estimation};
\coordinate (motorbend) at (170,-10.5);
\coordinate (motorup) at (170,9);
\draw[flow] (attctrl.east) -- (motorbend) -- (motorup) -- (plant.west);
\node[signal,align=center] at (165.5,-8.2) {motor\\commands};
\coordinate (sensorbranch) at (186,-1.5);
\draw[line width=0.65pt] (plant.south) -- (sensorbranch);
\fill (sensorbranch) circle (0.65pt);
\draw[flow] (sensorbranch) -| (imu.north);
\draw[flow] (sensorbranch) -| (fastlio.north);
\draw[flow] (imu.south) -- (imu.south |- ekf.north);
\draw[flow] (fastlio.south) -- (fastlio.south |- ekf.north);
\coordinate (attfbentry) at (144,-21.5);
\draw[line width=0.65pt] (ekf.west) --
  node[signal,align=center,pos=0.28,below=0.9mm] {attitude, rates} (attfbentry);
\draw[flow] (attfbentry) -- (attctrl.south);
\draw[flow,crossing]
  (ekf.south) -- ++(0,-2) -| node[signal,pos=0.46,above]
  {state feedback $X_0$} (ekftarget);

\end{tikzpicture}%
}
\caption{Controller architecture. The NMPC solver combines the exact
RK4-discretized model, FBDE derivatives, the OCP update, and a shifted warm
start.}
\label{fig:method_framework}
\end{figure*}

Fig.~\ref{fig:method_framework} shows how state estimation, the proposed
solver, and control execution are connected.

The remainder of this article is organized as follows.
The quadrotor system model, its exact discrete-time representation, and the
finite-horizon NMPC problem are presented in
Section~\ref{sec:formulation}. The FBDE-based gradient and Hessian recursions,
together with the OCP solver,
are developed in Section~\ref{sec:method}. Gazebo validation and real-flight
comparisons are reported in Sections~\ref{sec:results}
and~\ref{sec:flight}, respectively. Finally, conclusions are drawn in
Section~\ref{sec:conclusion}.

\medskip
\noindent\textbf{Notation:}
$\R^n$ denotes the $n$-dimensional Euclidean space.
The symbols $(\cdot)^\top$ and $\norm{\cdot}$ denote transposition and the
Euclidean norm, respectively, while $0_3$ and $I_3$ denote the
$3\times3$ zero and identity matrices.

\section{Problem Formulation}
\label{sec:formulation}

This section develops the discrete translational model used for prediction.
It then formulates the finite-horizon NMPC problem and eliminates the predicted
states to leave the input sequence as the decision variable.

\subsection{Model and Exact Discretization}

\begin{figure}[H]
\centering
\vspace*{3mm}
\resizebox{0.92\columnwidth}{!}{%
\begin{tikzpicture}[
  x=1mm,y=1mm,>=Latex,
  vector/.style={->,line width=1.15pt},
  frameaxis/.style={->,line width=1.0pt},
  axislabel/.style={font=\small\bfseries,fill=white,inner sep=1.0pt},
  varlabel/.style={font=\small\bfseries,fill=white,inner sep=1.1pt}
]
  \path[use as bounding box] (0,0) rectangle (100,66);
  \fill[white] (0,0) rectangle (100,66);
  \node[inner sep=0] (uavimage) at (51,33)
    {\includegraphics[width=50mm]{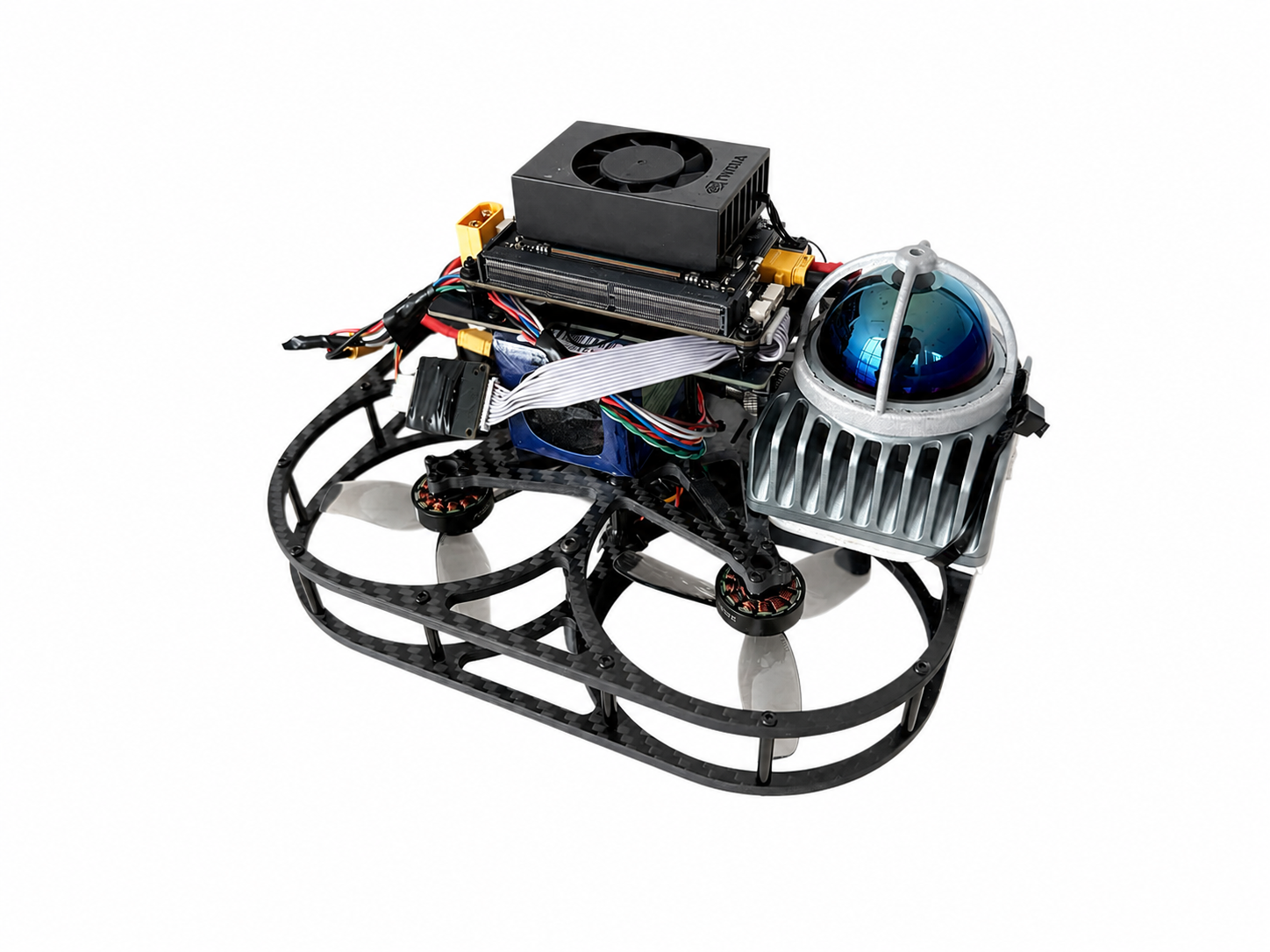}};

  \coordinate (OI) at (9,7);
  \draw[frameaxis,black] (OI) -- (28,7)
    node[axislabel,right] {$x_I$};
  \draw[frameaxis,black] (OI) -- (9,27)
    node[axislabel,above] {$z_I$};
  \draw[frameaxis,black] (OI) -- (1.5,0.8)
    node[axislabel,below left] {$y_I$};
  \node[varlabel,anchor=north west] at (9,6.5) {$O_I$};

  \coordinate (OB) at (52,36);
  \fill[white] (OB) circle (1.7pt);
  \draw[black,line width=0.7pt] (OB) circle (1.7pt);
  \draw[black!65,line width=0.7pt] (OB) -- (72,40)
    node[varlabel,right] {$O_b$};
  \coordinate (XBend) at (75,32);
  \coordinate (YBend) at (40,19);
  \draw[frameaxis,blue!80!black] (OB) -- (XBend)
    node[axislabel,right] {$x_b$};
  \draw[frameaxis,blue!80!black] (OB) -- (YBend)
    node[axislabel,below right] {$y_b$};
  \draw[frameaxis,blue!80!black] (OB) -- (52,52)
    node[axislabel,above] {$z_b$};

  \draw[vector,dashed,magenta!85!black] (OI) -- (OB);
  \node[varlabel,text=magenta!85!black] at (31,13) {$p=[x,y,z]^\top$};
  \draw[vector,green!55!black] (OB) -- (77,48)
    node[varlabel,above right] {$v=[v_x,v_y,v_z]^\top$};
  \draw[vector,red!80!black,line width=1.45pt] (53,55) -- (53,65)
    node[varlabel,above] {$T$};
  \coordinate (rollcenter) at ($(OB)!0.75!(XBend)$);
  \begin{scope}[rotate around={-9.87:(rollcenter)}]
    \draw[red!55!black,dashed,line width=0.75pt]
      ($(rollcenter)+(2,0)$)
      arc[start angle=0,end angle=180,x radius=2,y radius=4];
    \draw[blue!80!black,line width=1.0pt]
      ($(rollcenter)+(-3,0)$) -- ($(rollcenter)+(3,0)$);
    \draw[->,red!75!black,line width=1.05pt]
      ($(rollcenter)+(-2,0)$)
      arc[start angle=180,end angle=350,x radius=2,y radius=4];
  \end{scope}
  \node[varlabel,text=red!75!black] at (68,27) {$\phi$};

  \coordinate (pitchcenter) at ($(OB)!0.78!(YBend)$);
  \begin{scope}[rotate around={54.78:(pitchcenter)}]
    \draw[orange!60!black,dashed,line width=0.75pt]
      ($(pitchcenter)+(2,0)$)
      arc[start angle=0,end angle=180,x radius=2,y radius=4];
    \draw[->,orange!85!black,line width=1.05pt]
      ($(pitchcenter)+(-2,0)$)
      arc[start angle=180,end angle=350,x radius=2,y radius=4];
    \draw[blue!80!black,line width=1.0pt]
      ($(pitchcenter)+(-4,0)$) -- ($(pitchcenter)+(4,0)$);
  \end{scope}
  \node[varlabel,text=orange!85!black] at (48,20) {$\theta$};
\end{tikzpicture}%
}
\caption{Coordinate frames and variables of the experimental
quadrotor. $O_I$ and $O_b$ denote the origins of the inertial and body-fixed
frames,
respectively. The total thrust acts along $z_b$; $\phi$ and
$\theta$ denote rotations about $x_b$ and $y_b$, respectively.}
\label{fig:quadrotor_coordinates}
\end{figure}

The coordinate convention is shown in Fig.~\ref{fig:quadrotor_coordinates}, where
$p=[x,y,z]^\top$, $v=[v_x,v_y,v_z]^\top$, and $U=[T,\phi,\theta]^\top$ denote
the position, velocity, and control input, respectively. Here, $T$ is the total
thrust, while $\phi$ and $\theta$ are the commanded angles tracked by the faster
attitude loop. With yaw regulated independently, the dynamics are
\begin{equation}
\left\{
\begin{aligned}
\ddot{x}&=\tfrac{T}{m}\sin\theta\cos\phi,\\
\ddot{y}&=-\tfrac{T}{m}\sin\phi,\\
\ddot{z}&=\tfrac{T}{m}\cos\theta\cos\phi-g,
\end{aligned}
\right.
\label{eq:cont}
\end{equation}
where $m$ is the quadrotor mass and $g$ is the gravitational acceleration.

Defining the state and input vectors as
\[
X=[x,y,z,v_x,v_y,v_z]^\top\in\R^6,
U=[T,\phi,\theta]^\top\in\R^3,
\]
the dynamics in \eqref{eq:cont} can be written as
\begin{equation}
\dot X=f(X,U)=A_cX+b(U),
\label{eq:affine}
\end{equation}
where
\begin{equation*}
A_c=
\begin{bmatrix}
0_3&I_3\\
0_3&0_3
\end{bmatrix},
b(U)=
\begin{bmatrix}
0\\
0\\
0\\
\tfrac{T}{m}\sin\theta\cos\phi\\
-\tfrac{T}{m}\sin\phi\\
\tfrac{T}{m}\cos\theta\cos\phi-g
\end{bmatrix}.
\end{equation*}

The drift matrix in \eqref{eq:affine} satisfies $A_c^2=0$. Let $h>0$ be the
prediction interval. Holding $U_k$ constant over $[kh,(k+1)h)$ also makes
$b(U_k)$ constant, and the RK4 stages give
\begin{equation}
F(X_k,U_k)=
\left[\begin{array}{@{}l@{}}
x_k+h\,v_{x,k}+\tfrac{h^2T_k}{2m}\sin\theta_k\cos\phi_k\\[1pt]
y_k+h\,v_{y,k}-\tfrac{h^2T_k}{2m}\sin\phi_k\\[1pt]
z_k+h\,v_{z,k}+\tfrac{h^2T_k}{2m}\cos\theta_k\cos\phi_k-\tfrac{h^2}{2}g\\[1pt]
v_{x,k}+\tfrac{hT_k}{m}\sin\theta_k\cos\phi_k\\[1pt]
v_{y,k}-\tfrac{hT_k}{m}\sin\phi_k\\[1pt]
v_{z,k}+\tfrac{hT_k}{m}\cos\theta_k\cos\phi_k-h g
\end{array}\right].
\label{eq:disc}
\end{equation}
The same update can be written more compactly as
\begin{align}
p_{k+1}&=p_k+h v_k+\frac{h^2}{2}a(U_k),
\label{eq:position_update}\\
v_{k+1}&=v_k+h a(U_k),
\label{eq:velocity_update}
\end{align}
where the acceleration vector is
\begin{equation*}
\begin{aligned}
a(U_k)&=[a_{x,k},a_{y,k},a_{z,k}]^\top,\\
a_{x,k}&=\tfrac{T_k}{m}\sin\theta_k\cos\phi_k,\\
a_{y,k}&=-\tfrac{T_k}{m}\sin\phi_k,\\
a_{z,k}&=\tfrac{T_k}{m}\cos\theta_k\cos\phi_k-g.
\end{aligned}
\end{equation*}
Because $a(U_k)$ is constant over the interval, velocity changes linearly and
integration adds $h^2a(U_k)/2$ to the position update. Thus,
\eqref{eq:disc} is the exact discrete-time solution of \eqref{eq:affine} under
held inputs. Forward Euler gives the same velocity update but omits this
acceleration-induced displacement, causing an $O(h^2)$ one-step position
error.

\subsection{NMPC Formulation}

For trajectory tracking, let $(X_k^r,U_k^r)$ denote the reference sequences
over a horizon of $N$ stages. The reference input $U_k^r$ is obtained from the
reference acceleration using the standard differential-flatness relation for
quadrotors~\cite{mellinger2011}. The corresponding tracking errors are
$X_k^e=X_k-X_k^r$ and $U_k^e=U_k-U_k^r$. Collecting the input sequence as
$\mathcal U=[U_0^\top,\ldots,U_{N-1}^\top]^\top\in\R^{3N}$ gives the following
NMPC problem:
\begin{equation}
\left\{
\begin{aligned}
&\min_{\mathcal U}\quad J,\\
&J=\sum_{k=0}^{N-1}\Big[(X^e_k)^\top QX^e_k+(U^e_k)^\top R_uU^e_k\Big]+(X^e_N)^\top Q_fX^e_N,\\
&\text{s.t.}\quad X_{k+1}=F(X_k,U_k),\quad X_0\ \text{given},
\end{aligned}
\right.
\label{eq:ocp}
\end{equation}
where $Q\succeq0$, $Q_f\succeq0$, $R_u\succ0$, and $F$ is given by
\eqref{eq:disc}. Repeated application of the system model gives
\begin{equation}
X_k=\Psi_k(X_0,U_0,\dots,U_{k-1}),
\label{eq:state_transition_model}
\end{equation}
where $\Psi_k$ denotes the resulting $k$-step system model. Substituting \eqref{eq:state_transition_model} into the cost \eqref{eq:ocp} yields
\begin{equation}
\min_{\mathcal U\in\R^{3N}}J(X_0,\mathcal U).
\label{eq:op}
\end{equation}
Eq.~\eqref{eq:op} is the problem solved at each sampling instant. Its only
decision variable is the physical input sequence $\mathcal U$; the state
evolution remains embedded in $\Psi_k$. This form makes it possible to update
the complete sequence using derivatives of $J$ with respect to $\mathcal U$.

\section{Solver and Analytical Derivatives}
\label{sec:method}

To solve~\eqref{eq:op}, the OCP iteration requires the gradient and Hessian of
$J(X_0,\mathcal U)$. The update is stated first, followed by compact FBDE expressions, the shifted warm start, and the complete online algorithm.

\subsection{OCP Iterative Update}
At optimization iteration $i$, the current input sequence $\mathcal U_i$ is
updated by
\begin{equation}
\mathcal U_{i+1}=\mathcal U_i-d_i(\mathcal U_i),
\label{eq:general_iterative_update}
\end{equation}
where $d_i(\mathcal U_i)$ denotes the update vector computed from the
derivatives of $J(X_0,\mathcal U)$.

The OCP method in~\cite{zhang2024scis,lv2026tie} constructs this update
recursively using the gradient and Hessian.
For the present NMPC problem, the cost, gradient, and Hessian at
$\mathcal U_i$ are denoted by
\begin{equation*}
\begin{aligned}
J_i&=J(X_0,\mathcal U_i),\\
\nabla J_i&=\nabla_{\mathcal U}J(X_0,\mathcal U_i),\\
H_i&=\nabla_{\mathcal U}^2J(X_0,\mathcal U_i).
\end{aligned}
\end{equation*}
For a fixed matrix $\mathcal{R}\succ0$, the search direction is
generated by
\begin{align}
d_0(\mathcal U_i)&=(\mathcal{R}+H_i)^{-1}\nabla J_i,
\label{eq:direction_initial}\\
d_j(\mathcal U_i)&=(\mathcal{R}+H_i)^{-1}
 \big[\nabla J_i+\mathcal{R}d_{j-1}(\mathcal U_i)\big],
 j=1,\ldots,i,
\label{eq:direction_recursion}\\
\mathcal U_{i+1}&=\mathcal U_i-d_i(\mathcal U_i).
\label{eq:unconstrained_update}
\end{align}
Eq.~\eqref{eq:direction_initial} gives the initial direction. It is refined
$i$ times by~\eqref{eq:direction_recursion}, after which
\eqref{eq:unconstrained_update} advances the input sequence. Thus the inner
recursion grows with the outer iteration index. Since $\mathcal R+H_i$ remains
fixed during these refinements, its inverse is computed once per outer
iteration and reused.

The required $\nabla J_i$ and $H_i$ are derived next. The gradient follows
from the state and costate equations; differentiating these recursions gives
the Hessian.

\subsection{Analytical Gradient}
After forward prediction with~\eqref{eq:disc}, define
$E_k=\partial F/\partial U_k$ and
$A_d=\partial F/\partial X=\left[\begin{smallmatrix}I_3&hI_3\\0_3&I_3\end{smallmatrix}\right]$.
The gradient blocks and backward costate recursion are
\begin{equation*}
\begin{aligned}
\gamma_k&=2R_uU_k^e+E_k^\top\lambda_{k+1},\\
\lambda_N&=2Q_fX_N^e,\\
\lambda_k&=2QX_k^e+A_d^\top\lambda_{k+1},
\quad k=N-1,\ldots,0.
\end{aligned}
\end{equation*}
The complete gradient is
$\nabla_{\mathcal U}J=[\gamma_0^\top,\ldots,\gamma_{N-1}^\top]^\top$.
For $Q=\operatorname{diag}(q_1,\ldots,q_6)$, let
$Q_p=\operatorname{diag}(q_1,q_2,q_3)$ and
$Q_v=\operatorname{diag}(q_4,q_5,q_6)$. Partitioning the costate into
position and velocity components gives
\begin{equation}
\begin{aligned}
\lambda_{p,k}&=2Q_pp_k^e+\lambda_{p,k+1},\\
\lambda_{v,k}&=2Q_vv_k^e+h\lambda_{p,k+1}+\lambda_{v,k+1}.
\end{aligned}
\label{eq:costate}
\end{equation}

The remaining term follows directly from the position and velocity updates:
\begin{equation*}
E_k^\top\lambda_{k+1}
=h\left(\frac{\partial a(U_k)}{\partial U_k}\right)^\top
\left(\lambda_{v,k+1}+\frac h2\lambda_{p,k+1}\right).
\end{equation*}
For compactness, define
\begin{equation*}
\mu_{k+1}
=\lambda_{v,k+1}+\frac h2\lambda_{p,k+1}
=\begin{bmatrix}\mu_{4,k+1}&\mu_{5,k+1}&\mu_{6,k+1}\end{bmatrix}^\top.
\end{equation*}
Thus,
\begin{equation*}
\gamma_k=2R_uU_k^e
+h\left(\frac{\partial a(U_k)}{\partial U_k}\right)^\top\mu_{k+1}.
\end{equation*}
For $R_u=\operatorname{diag}(R_{u,1},R_{u,2},R_{u,3})$, expanding the input
derivative gives
\begin{align}
\gamma_{T,k}&=2R_{u,1}U^e_{1,k}
+\tfrac{h}{m}(\sin\theta_k\cos\phi_k\,\mu_{4,k+1}\notag\\
&\quad-\sin\phi_k\,\mu_{5,k+1}
+\cos\theta_k\cos\phi_k\,\mu_{6,k+1}),\notag\\
\gamma_{\phi,k}&=2R_{u,2}U^e_{2,k}
+\tfrac{hT_k}{m}(-\sin\theta_k\sin\phi_k\,\mu_{4,k+1}\notag\\
&\quad-\cos\phi_k\,\mu_{5,k+1}
-\cos\theta_k\sin\phi_k\,\mu_{6,k+1}),\notag\\
\gamma_{\theta,k}&=2R_{u,3}U^e_{3,k}
+\tfrac{hT_k}{m}(\cos\theta_k\cos\phi_k\,\mu_{4,k+1}\notag\\
&\quad-\sin\theta_k\cos\phi_k\,\mu_{6,k+1}).
\label{eq:gradient_components}
\end{align}
Here, $\gamma_k=[\gamma_{T,k},\gamma_{\phi,k},\gamma_{\theta,k}]^\top$.
Thus, one forward state evaluation and one backward costate evaluation give
the complete gradient, without explicit differentiation of $\Psi_k$.

\subsection{Exact Hessian}
The Hessian follows by differentiating the state and costate recursions.
For input stage $r$, let
$S_j^{(r)}=\partial X_j/\partial U_r$ and
$P_j^{(r)}=\partial\lambda_j/\partial U_r$.
Since $A_d$ is constant, these $6\times3$ sensitivities satisfy
\begin{equation*}
\begin{aligned}
S_{j+1}^{(r)}&=A_dS_j^{(r)}+\delta_{jr}E_r,
& S_0^{(r)}&=0,\\
P_j^{(r)}&=2QS_j^{(r)}+A_d^\top P_{j+1}^{(r)},
& P_N^{(r)}&=2Q_fS_N^{(r)},
\end{aligned}
\end{equation*}
Here $j=0,\ldots,N-1$ and $\delta_{jr}$ is the Kronecker delta. The first recursion is evaluated
forward and the second backward. Differentiating $\gamma_k$ gives the
$3\times3$ Hessian blocks
\begin{equation}
H_{kr}=E_k^\top P_{k+1}^{(r)}+\delta_{kr}(2R_u+\Theta_k),
\label{eq:hess}
\end{equation}
with the local curvature term
\begin{equation}
\Theta_k=h\sum_{a=1}^{3}\mu_{a+3,k+1}
\frac{\partial^2 a_a(U_k)}{\partial U_k^2}.
\label{eq:theta}
\end{equation}
Here $a_a$ is component $a$ of the acceleration vector. The derivatives
in~\eqref{eq:theta} are evaluated analytically from the trigonometric
expressions in Section~\ref{sec:formulation}.
Because the model has no state--state or state--input second derivatives,
\eqref{eq:hess} includes all second-order terms; retaining $\Theta_k$
distinguishes the exact Hessian from a Gauss--Newton approximation.

For this translational model, the future-state terms in the sensitivity
recursions can also be summed in closed form. Define
\begin{equation*}
c_{a,k}=\frac{\partial[h a_a(U_k)]}{\partial U_k}
\in\mathbb R^{1\times3},\quad a=1,2,3.
\end{equation*}
Let $Q_f=\operatorname{diag}(q_{f,1},\ldots,q_{f,6})$. For
$0\le r\le s\le N-1$, set $M=N-s-1$ and $\eta=s-r$. The $M$ shared stage
states are indexed by $\nu=0,\ldots,M-1$, and the terminal state corresponds
to $\nu=M$. Define
\begin{align*}
A_{M,\eta}&=M(4M^2+6\eta M-1)/12,\\[-2pt]
B_{M,\eta}&=(M+\eta+1/2)(M+1/2).
\end{align*}
Direct evaluation of the shared future-state contributions gives
\begin{align}
w_{a,rs}
&=2q_a h^2A_{M,\eta}+2q_{a+3}M\notag\\
&\quad+2q_{f,a}h^2B_{M,\eta}+2q_{f,a+3},
\quad a=1,2,3.
\label{eq:closed_form_weight}
\end{align}
For $r\le s$, the corresponding exact Hessian blocks are
\begin{equation}
\begin{aligned}
H_{rs}&=\sum_{a=1}^{3}w_{a,rs}c_{a,r}^\top c_{a,s}
 +\begin{cases}0,&r<s,\\2R_u+\Theta_r,&r=s,\end{cases}\\
H_{sr}&=H_{rs}^\top.
\end{aligned}
\label{eq:scalar_exact_hessian}
\end{equation}
Eqs.~\eqref{eq:theta} and~\eqref{eq:scalar_exact_hessian} then assemble the
exact Hessian from scalar terms after the state and costate evaluations,
without explicitly differentiating $\Psi_k$.

\subsection{Shifted Warm Start}
At sampling instant $t$, discard the previously applied input, shift the
remaining sequence by one stage, and append the current terminal reference:
\begin{equation}
\begin{aligned}
U_{k|t}^{0}&=
\begin{cases}
U_{k+1|t-1}^{\star},&k=0,\ldots,N-2,\\
U^r_{N-1|t},&k=N-1,
\end{cases}\\
\mathcal U_t^0&=
\begin{bmatrix}
(U_{0|t}^{0})^\top&\cdots&(U_{N-1|t}^{0})^\top
\end{bmatrix}^\top.
\end{aligned}
\label{eq:shifted_warm_start}
\end{equation}
In the present controller, $U^r_{N-1|t}$ is generated from the terminal
reference acceleration by differential flatness. It is used only to complete
the initial guess; after optimization, only the first optimized input
$U_{0|t}^{\star}$ is applied. Since two successive NMPC problems differ only by
one sampling interval, the shifted sequence is usually close to the new
solution and reduces the number of online iterations.

\subsection{Overall Algorithm}
Algorithm~\ref{alg:1} collects the operations at one sampling instant. From
the shifted warm start, the state and costate equations give $\nabla J_i$;
the closed-form Hessian expressions give $H_i$; and the nested direction recursion
updates the input sequence.

\begin{algorithm}[t]
\caption{Numerical Algorithm for Solving the NMPC Problem}
\label{alg:1}
\begin{algorithmic}[1]
\State \textbf{Input:} Measured state $X_0$; reference sequences $X^r,U^r$;
       previous solution $\mathcal U_{t-1}^{\star}$;
       maximum outer iterations $i_{\max}$; gradient stopping tolerance
       $\varepsilon_g$; positive-definite regularization matrix
       $\mathcal{R}\succ0$
\State \textbf{Output:} Final input sequence $\mathcal U_t^\star$ and first
       command $U_{0|t}^\star=(T_t,\phi_{d,t},\theta_{d,t})$
\State Discard the applied $U_{0|t-1}^{\star}$
\State Set $U_{k|t}^{0}\gets U_{k+1|t-1}^{\star}$ for $k=0,\ldots,N-2$
\State Set $U_{N-1|t}^{0}\gets U^r_{N-1|t}$ and stack
       $\mathcal U_0=\mathcal U_t^0$ according to~\eqref{eq:shifted_warm_start}
\State $\mathcal U_t^\star\gets\mathcal U_0$
\For{$i=0$ to $i_{\max}-1$}
  \State Evaluate $\nabla J_i=\nabla_{\mathcal U}J(X_0,\mathcal U_i)$ using
         \eqref{eq:disc}, \eqref{eq:costate}, and
         \eqref{eq:gradient_components}
  \If{$\norm{\nabla J_i}_\infty<\varepsilon_g$}
    \State $\mathcal U_t^\star\gets\mathcal U_i$; \textbf{break}
  \EndIf
  \State Compute $H_i=\nabla_{\mathcal U}^2J(X_0,\mathcal U_i)$ from
         \eqref{eq:theta} and \eqref{eq:scalar_exact_hessian}
  \State $d_0(\mathcal U_i)=(\mathcal{R}+H_i)^{-1}\nabla J_i$
  \If{$i\geq1$}
    \For{$j=1$ to $i$}
      \State $d_j(\mathcal U_i)=(\mathcal{R}+H_i)^{-1}
             [\nabla J_i+\mathcal{R}d_{j-1}(\mathcal U_i)]$
    \EndFor
  \EndIf
  \State $\mathcal U_{i+1}=\mathcal U_i-d_i(\mathcal U_i)$
  \State $\mathcal U_t^\star\gets\mathcal U_{i+1}$
\EndFor
\State Apply $U_{0|t}^\star$ and retain $\mathcal U_t^\star$ for the next
       sampling instant
\end{algorithmic}
\end{algorithm}

\medskip
\noindent\textbf{Remark 1 (Extension to Input Constraints):}
Algorithm~\ref{alg:1} solves the unconstrained NMPC problem. If box constraints
$\underline{\mathcal U}\leq\mathcal U\leq\overline{\mathcal U}$ are required,
the same analytical gradient and Hessian can be combined with a primal
active-set method. At iteration $i$, components at the lower bound with
$[\nabla J_i]_r>0$ and those at the upper bound with $[\nabla J_i]_r<0$ are
placed in the active set; all remaining components form the free set. The
active components are fixed, and the nested direction recursion is applied to
the gradient and Hessian restricted to the free variables. A feasible step is
taken until a new bound is reached, while an active component is released when
its gradient no longer satisfies the corresponding sign condition. Thus,
constraint handling changes only the active-set identification and the reduced
direction calculation without affecting derivative computation through the FBDE recursions.

\section{Gazebo Validation}
\label{sec:results}

\subsection{Unified Setup and Compared Solvers}
The evaluation proceeds from derivative checks to controlled discretization
and warm-start ablations, followed by a closed-loop timing comparison at a
similar tracking-accuracy level.
The Gazebo quadrotor tracks the three-dimensional Lissajous reference
$x^r=3+2\sin(\omega t)$, $y^r=3+2\sin(2\omega t)$, and
$z^r=5+1.5\sin(3\omega t)$, where $\omega=0.15$. The reference input $U^r$ is
obtained from the reference acceleration by differential flatness. All controllers operate at $50$\,Hz with $h=0.02$\,s and $N=20$.
In both the Gazebo and real-flight tests, Algorithm~\ref{alg:1} uses
$i_{\max}=100$, $\varepsilon_g=10^{-6}$, and
$\mathcal R=0.01I_{3N}$.

The comparison includes Algorithm~\ref{alg:1}, iLQR, acados SQP+HPIPM, and
acados SQP+qpOASES. All four solve the same unconstrained NMPC problem and use
the same Gazebo plant, reference, prediction horizon, inner-loop controller,
thrust-command mapping factor of $0.785$, and ROS~2 interfaces. The common cost
weights are
$Q=Q_f=\operatorname{diag}(10.5,10.5,12,3.8,3.8,5)$ and
$R_u=\operatorname{diag}(0.08,8,12)$. The corresponding NMPC cost is given
in~\eqref{eq:ocp}. Solver-specific numerical settings were
tuned separately for good closed-loop performance, while all methods use the
same stopping condition $\norm{\nabla_{\mathcal U}J}_\infty<10^{-6}$.
Each method is evaluated in ten independent closed-loop trials. The reported
trajectories and solver-time histories are averaged over the ten trials, while
each tabulated statistic is first computed for an individual trial and then
averaged over the ten trials. For all four methods, the solver execution time is
measured within the C implementation, from entry to exit of the numerical
solver, excluding reference generation, state estimation, ROS~2 communication,
and command publication. This measurement is hereafter referred to as the
C-core solution time.

A central-difference check gives relative errors of $5.286\times10^{-9}$ and
$1.440\times10^{-10}$ for the gradient and complete Hessian, respectively. The
Hessian is symmetric to the reported numerical precision, confirming the FBDE
recursions and the closed-form weights in~\eqref{eq:closed_form_weight}.

\subsection{Discretization and Warm-Start Ablations}
For the ablations, the same $2828$ post-transient states and references are
replayed with different discretizations and initializations. This removes
run-to-run plant variation from the comparison. The replay is offline; its
times are not treated as closed-loop solver measurements.

\begin{figure}[t]
\centering
\vspace*{1mm}
\begin{minipage}[t]{0.485\columnwidth}
  \centering
  \includegraphics[width=\linewidth,height=3.5cm]{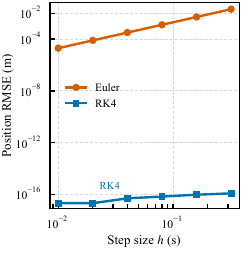}\\[-1mm]
  \footnotesize (a)
\end{minipage}\hfill
\begin{minipage}[t]{0.485\columnwidth}
  \centering
  \includegraphics[width=\linewidth,height=3.5cm]{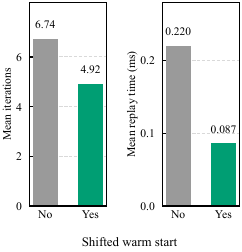}\\[-1mm]
  \footnotesize (b)
\end{minipage}
\caption{Controlled same-state replay ablations using states recorded from the
Gazebo closed-loop experiment for Algorithm~\ref{alg:1}:
(a) one-step position error of the RK4 and forward-Euler models; (b) mean
iteration count and offline-replay solver time with cold and shifted warm
starts.}
\label{fig:gazebo_ablations}
\end{figure}

As shown in Fig.~\ref{fig:gazebo_ablations}(a), the one-step position error of
forward Euler scales as $O(h^2)$ with the step size $h$. The RK4 error remains
at floating-point level because the
held-input dynamics are integrated exactly by the RK4 model.
Fig.~\ref{fig:gazebo_ablations}(b) shows that the shifted warm start reduces
the mean iteration count from $6.738$ to $4.925$ and the mean replay time from
$0.2202$ to $0.08656$~ms. These changes correspond to reductions of $26.91\%$
and $60.70\%$, respectively. Each solve is repeated five times with alternating
call order, and reference generation is excluded from the timing. The absolute
replay times in Fig.~\ref{fig:gazebo_ablations} are therefore not compared
directly with Table~\ref{tab:gazebo_all}.

\subsection{Solver-Time and Iteration Comparison}
\begin{figure}[t]
\centering
\includegraphics[width=\columnwidth,height=6.5cm]{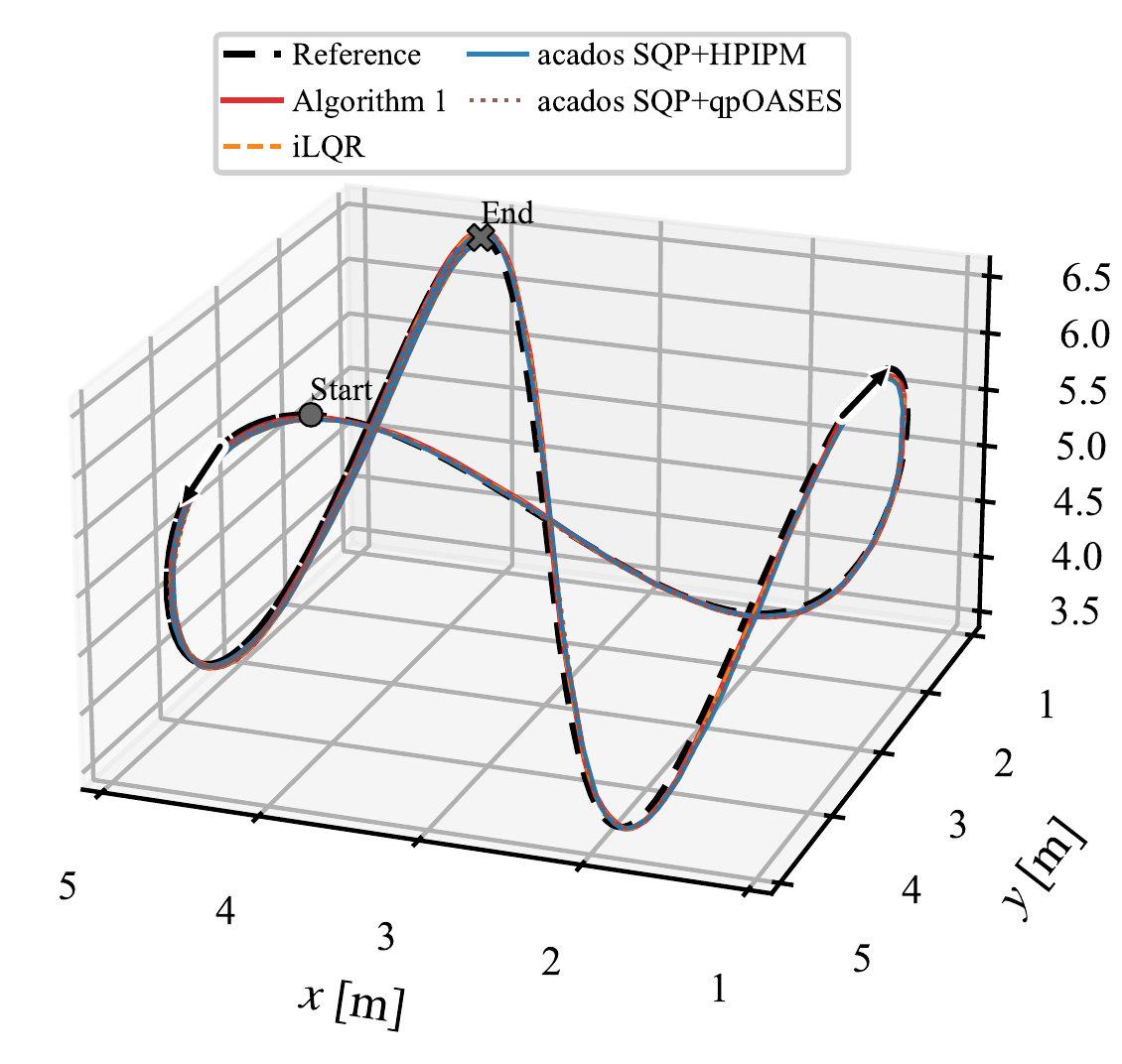}
\caption{Gazebo three-dimensional trajectory tracking for the four methods.}
\label{fig:gazebo_3d}
\end{figure}

\begin{figure}[t]
\centering
\includegraphics[width=\columnwidth,height=4.5cm]{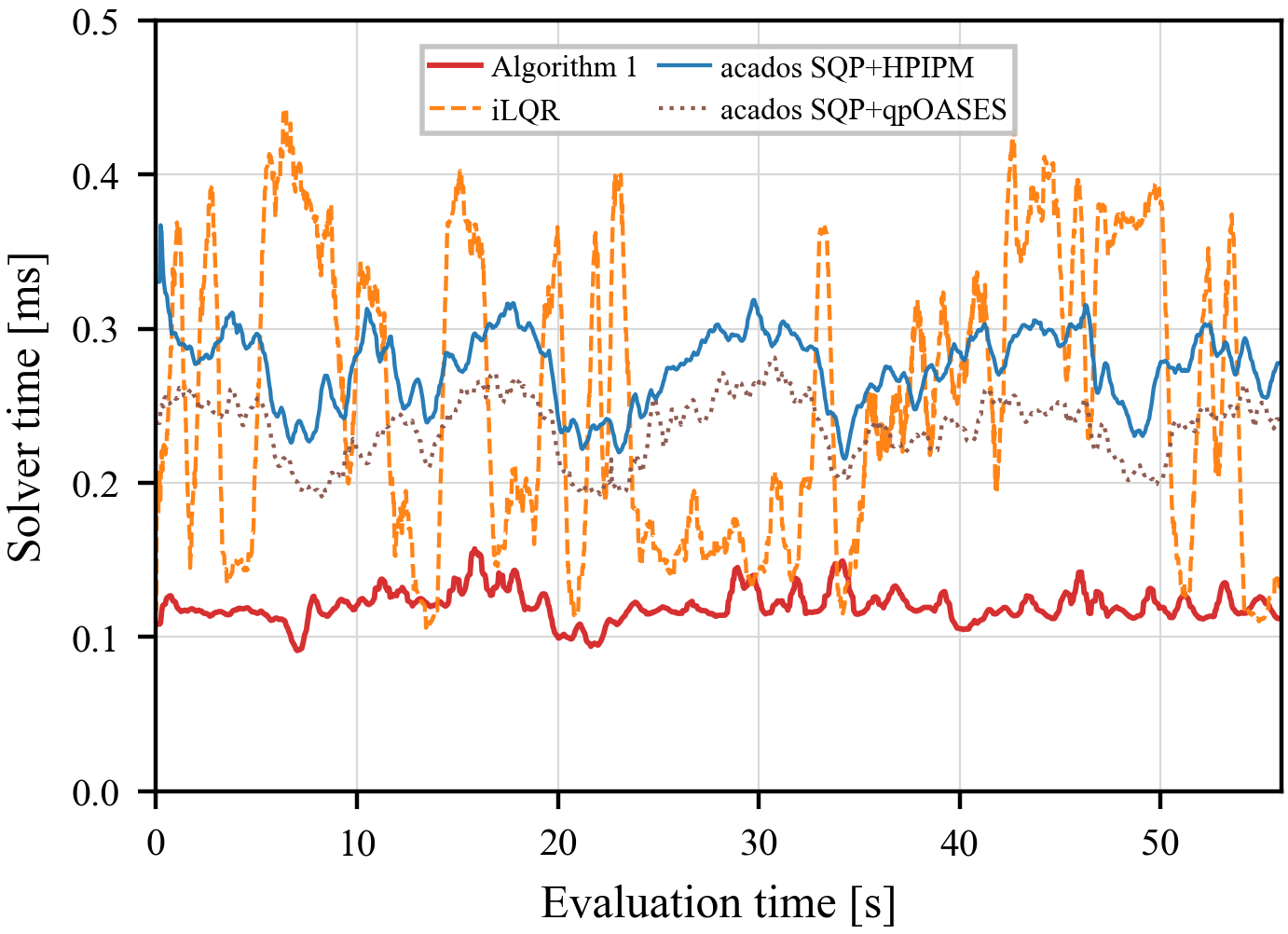}
\caption{Gazebo C-core solver-time histories for the four methods.}
\label{fig:gazebo_time}
\end{figure}

All four methods complete the three-dimensional trajectory in
Fig.~\ref{fig:gazebo_3d}. Their 3-D RMSE values differ by at most
$4.857$\,mm, providing a comparable tracking level for the solver-time
comparison.
With this tracking comparison established, Fig.~\ref{fig:gazebo_time} shows
how the C-core solution times vary during closed-loop operation. Algorithm~1
remains near $0.12$\,ms for most of the evaluated interval, below the other
three methods, while iLQR varies more widely. The unsmoothed statistics are
reported in Table~\ref{tab:gazebo_all}.

\begin{table}[t]
\centering
\caption{Gazebo Performance
(C-Core Times in ms)}
\label{tab:gazebo_all}
\scriptsize
\setlength{\tabcolsep}{2.0pt}
\renewcommand{\arraystretch}{1.02}
\begin{tabular*}{\columnwidth}{@{\extracolsep{\fill}}lcccc@{}}
\toprule
Method & Mean & P99 & Max. & \shortstack{Iterations\\mean/max.}\\
\midrule
\textbf{Algorithm 1}  & \textbf{0.1197}
& \textbf{0.2771} & \textbf{0.6049} & \textbf{4.92/5}\\
iLQR  & 0.2537 & 0.6659 & 0.9405 & 6.43/10\\
acados SQP+HPIPM  & 0.2748 & 0.3700 & 0.6940 & 2.66/3\\
acados SQP+qpOASES  & 0.2358 & 0.3540 & 0.5680 & 2.66/3\\
\bottomrule
\end{tabular*}
\end{table}

Algorithm~1 requires $4.92$ iterations on average and never more than five.
Its mean, 99th percentile (P99), and maximum solve times are $0.1197$, $0.2771$, and
$0.6049$\,ms, respectively. Relative to iLQR, acados SQP+HPIPM, and acados
SQP+qpOASES, its mean time is lower by $52.82\%$, $56.44\%$, and $49.24\%$,
or factors of $2.12$, $2.30$, and $1.97$. Iteration counts have different
meanings across the solvers, but Algorithm~1 meets its stopping condition
within five iterations at every sampling instant.

The timing difference is consistent with the solver design: analytical
derivatives avoid repeated numerical differentiation, the OCP refinements
reuse one inverse per outer iteration, and the shifted previous solution
initializes the next solve. No QP subproblem or general-purpose nonlinear
programming package is used inside the online solve. Whether this advantage
persists on the onboard processor is examined next.

\section{Real-Flight Validation}
\label{sec:flight}

\subsection{Flight Platform and Comparison Setup}
The real-flight tests examine whether the Gazebo timing advantage persists on
the onboard processor. The four methods in Table~\ref{tab:flight_timing}---Algorithm~1, iLQR,
acados SQP+HPIPM, and acados SQP+qpOASES---are tested onboard in ten
independent flights each. All four solve the same unconstrained NMPC problem.
Each controller tracks a circle of $1.0$\,m radius at an altitude
of approximately $0.5$\,m. The trajectory and controller parameters keep
the commanded thrust and attitude within the platform's safe operating
range. Each flight covers three $20$-s laps, with controller, reference, and
state-estimation data recorded at approximately $50$\,Hz. Reported
trajectories and timing statistics are averaged across the ten flights.

\begin{figure}[t]
\centering
\vspace*{1mm}
\includegraphics[width=5.5cm,height=3.5cm]{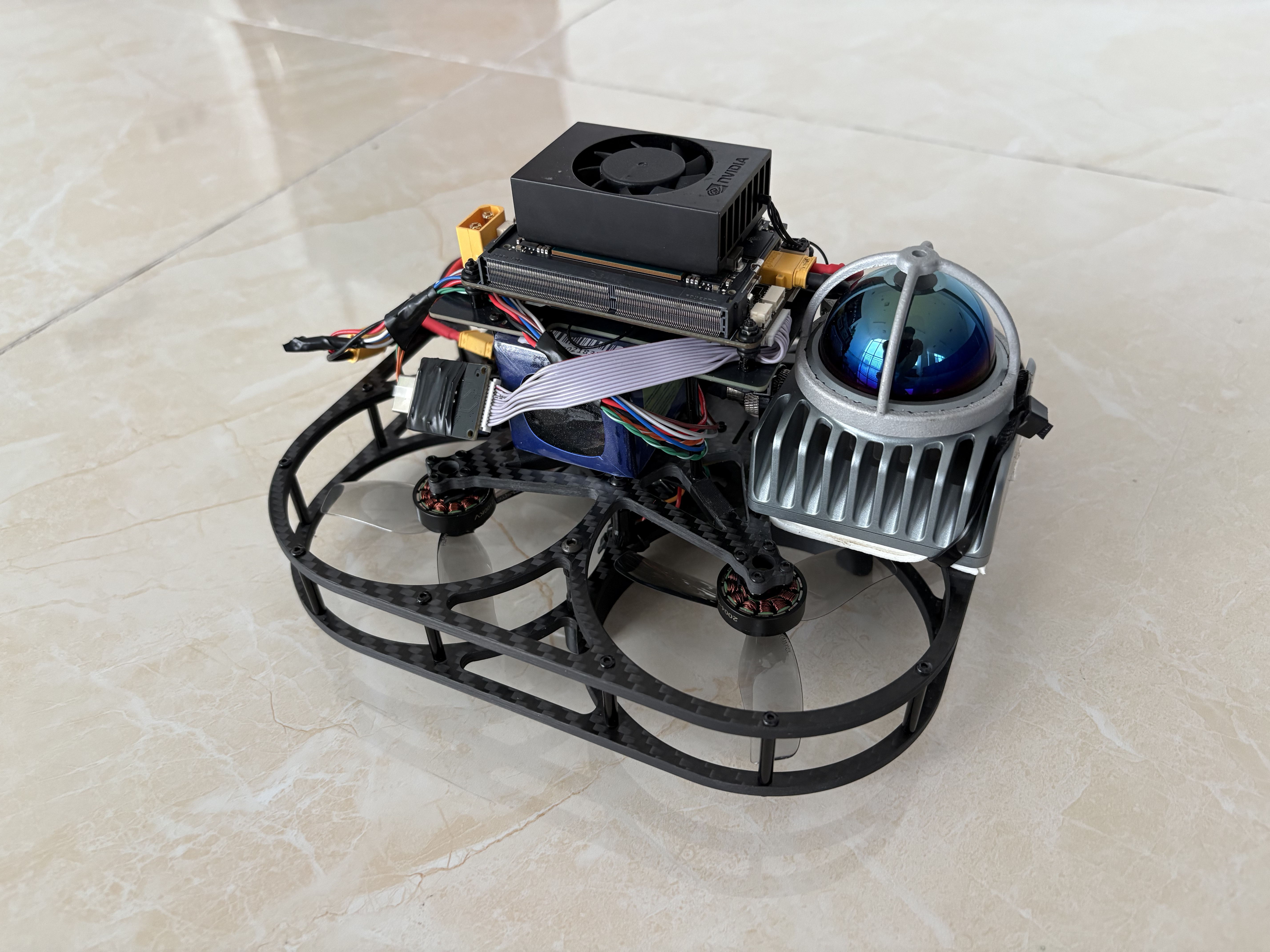}
\caption{Experimental quadrotor platform.}
\label{fig:quadrotor_platform}
\end{figure}

The experimental platform is shown in Fig.~\ref{fig:quadrotor_platform}. It is
built on a $3.5$-inch frame, with a MicoAir743v2-AIO-35A board providing flight
control and motor drive. All optimization methods run on the same NVIDIA
Jetson Orin NX $8$~GB computer through identical ROS~2 interfaces. LiDAR and
inertial measurements are processed by FAST-LIO to provide odometry, which is
fused with onboard IMU measurements by the EKF for state feedback. The airframe, sensors,
computer, communication interfaces, and inner-loop controller remain unchanged
across the runs; only the NMPC solution method is replaced.


\subsection{Solver-Time Comparison}
\begin{figure}[t]
\centering
\includegraphics[width=\columnwidth,height=6.5cm]{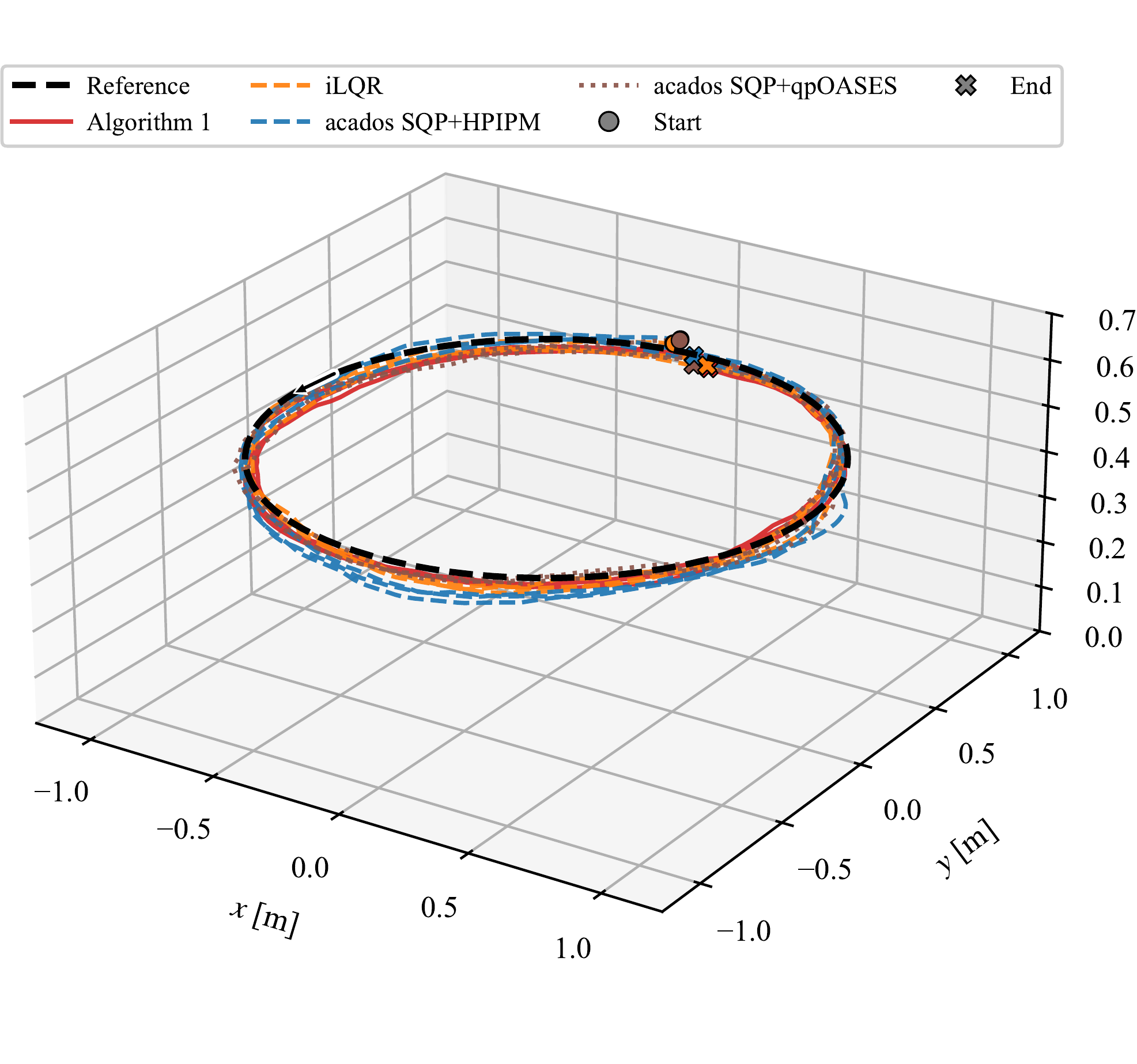}
\caption{Measured three-dimensional trajectory tracking for the four methods.}
\label{fig:acados_flight}
\end{figure}

In Fig.~\ref{fig:acados_flight}, all four controllers complete the circular
task in closed loop. Their 3-D tracking RMSE values are at most
$0.1174$\,m; the comparison below focuses on solver time within this
tracking-error range.

\begin{figure}[t]
\vspace*{1mm}
\centering
\includegraphics[width=\columnwidth,height=4.5cm]{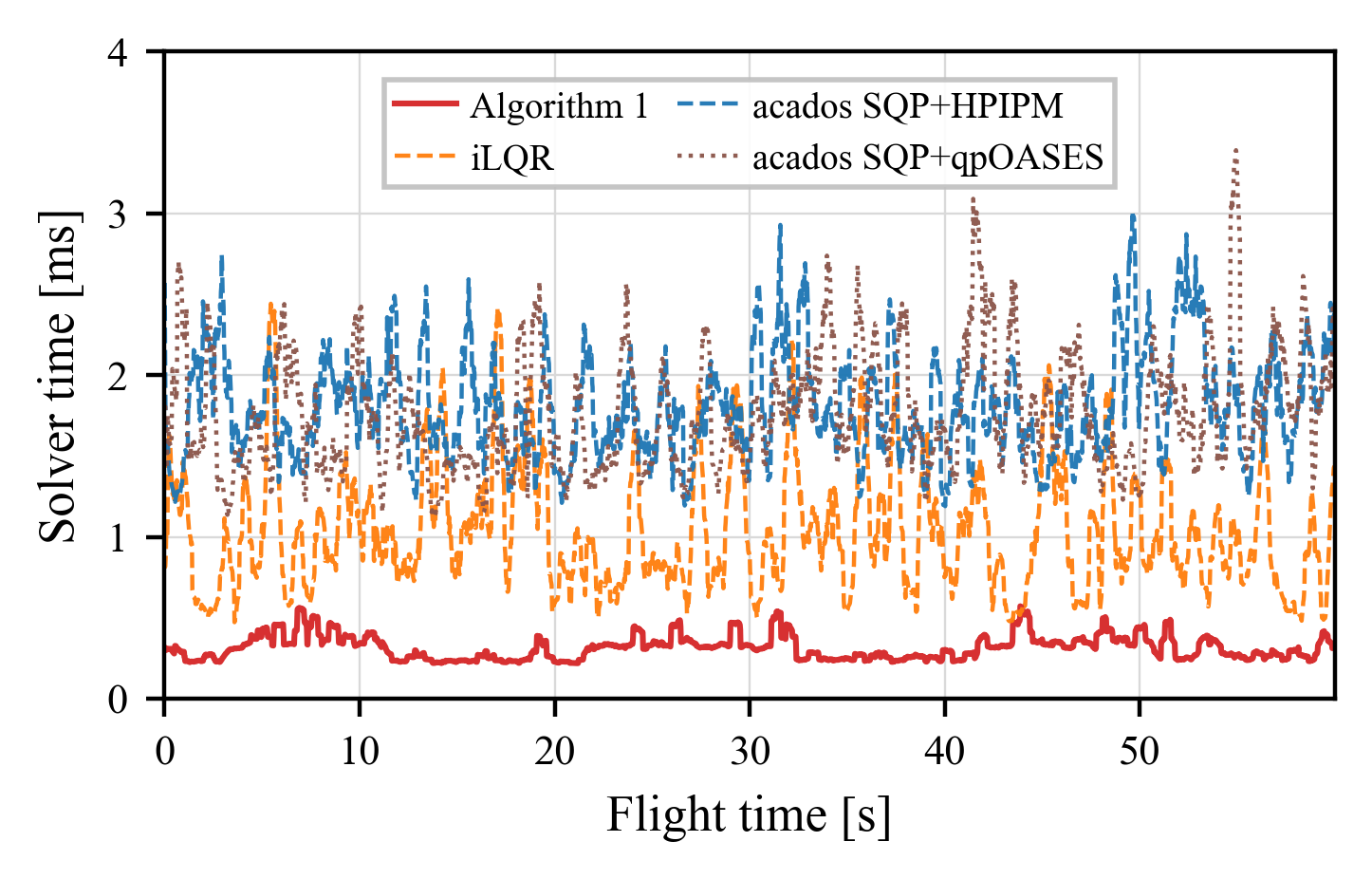}
\caption{Measured real-flight solver-time histories for the four methods.}
\label{fig:flight_solver_time}
\end{figure}

Fig.~\ref{fig:flight_solver_time} shows the measured C-core solve times over
the flight. Algorithm~1 is below the three baselines for most of the interval;
the difference is not confined to a few samples. Table~\ref{tab:flight_timing}
reports the unsmoothed statistics and iteration counts.

\begin{table}[t]
\centering
\caption{Measured Real-Flight Performance
(C-Core Times in ms)}
\label{tab:flight_timing}
\footnotesize
\setlength{\tabcolsep}{2pt}
\renewcommand{\arraystretch}{1.08}
\begin{tabular*}{\columnwidth}{@{\extracolsep{\fill}}lccccc@{}}
\toprule
Method & Mean & P95 & P99 & Max. & \shortstack{Iter.\\(mean/max.)}\\
\midrule
\textbf{Algorithm 1} & \textbf{0.3121} & \textbf{0.4316}
& \textbf{1.3105} & \textbf{3.7754} & \textbf{4.5/5}\\
iLQR & 1.0523 & 3.0614 & 5.5872 & 13.1394 & 4.7/11\\
acados SQP+HPIPM & 1.8014 & 4.5029 & 8.6506 & 14.9986 & 2.4/5\\
acados SQP+qpOASES & 1.7327 & 3.9870 & 7.5551 & 17.4178 & 2.4/5\\
\bottomrule
\end{tabular*}
\end{table}

Algorithm~1 takes $4.5$ iterations on average and at most five. Its mean,
95th percentile (P95), P99, and maximum solve times of $0.3121$, $0.4316$, $1.3105$, and
$3.7754$\,ms are the smallest in Table~\ref{tab:flight_timing}. Mean-time
speedup factors relative to iLQR, acados SQP+HPIPM, and acados SQP+qpOASES
are $3.37$, $5.77$, and $5.55$, respectively. The P99 times are lower by
$76.54\%$--$84.85\%$. The two acados configurations take fewer iterations on
average, yet require more time per completed NMPC solve.

The maximum time of each method remains below the $20$ ms sampling period. In
particular, the $3.7754$ ms worst-case time of Algorithm~1 occupies only
$18.9\%$ of this period and corresponds to a solver-only rate of approximately
$264.9$\,Hz. This value is not an end-to-end closed-loop rate because state
estimation, communication, and command publication also require time. The
measured mean and worst-case results leave substantial computational margin for
high-frequency onboard implementation.

The mean-time ordering is the same in Gazebo and flight: Algorithm~1 is
fastest, by factors of $1.97$--$2.30$ and $3.37$--$5.77$, respectively.
Different trajectories and operating environments preclude a direct comparison
of absolute times, but the advantage over the baselines persists onboard.

\section{Conclusion}\label{sec:conclusion}
This work addressed the online solution cost of UAV NMPC while retaining the
nonlinear translational predictor. The structure of that model permits an exact
discrete-time RK4 update and compact FBDE recursions for the gradient and exact
Hessian of the finite-horizon cost with respect to the input sequence. These derivatives feed the OCP nested direction recursion, with a shifted warm start carrying the previous solution into each new sampling instant. In Gazebo, the solver achieved lower C-core time than the
three baselines at similar tracking accuracy. The same timing advantage was
observed onboard: the mean solution-time speedups ranged from $3.37$ to $5.77$,
and the longest measured solve remained below the $20$ ms control period.
Future work will address input constraints and measure the complete
estimation--optimization--actuation latency at higher closed-loop rates.

\bibliographystyle{IEEEtran}
\bibliography{references}

\end{document}